\documentclass[10pt]{amsart}

\usepackage[matrix,arrow,curve,frame]{xy}
\usepackage{amsmath,amsthm,amssymb,enumerate}
\usepackage{latexsym}
\usepackage{amscd}

\usepackage[colorlinks=true]{hyperref}

\newtheorem{thm}[subsection]{Theorem}

\newtheorem{cor}[subsection]{Corollary}
\newtheorem{lemma}[subsection]{Lemma}

\newtheorem{remark}[subsection]{Remark}

\theoremstyle{definition}
\newtheorem{example}[subsection]{Example}

\newcommand{\cat}{\mathcal}
\newcommand{\lra}{\longrightarrow}
\newcommand{\lla}{\longleftarrow}
\newcommand{\llra}[1]{\stackrel{#1}{\lra}}
\newcommand{\llla}[1]{\stackrel{#1}{\lla}}

\newcommand{\R}{\mathbb R}
\newcommand{\Q}{\mathbb Q}
\newcommand{\Z}{\mathbb Z}
\newcommand{\N}{\mathbb N}
\newcommand{\C}{\mathbb C}
\newcommand{\CP}{\mathbb P}

\newcommand{\Gg}{{\cat G}}
\newcommand{\Jj}{{\cat J}}
\newcommand{\Oo}{{\cat O}}

\newcommand{\Om}{\Omega}
\newcommand{\om}{\omega}
\newcommand{\la}{\lambda}

\DeclareMathOperator{\Aut}{Aut}

\DeclareMathOperator{\Diff}{Diff}

\DeclareMathOperator{\Symp}{Symp}

\DeclareMathOperator{\Hol}{Hol}
\DeclareMathOperator{\Iso}{Iso}
\DeclareMathOperator{\Real}{Re}
\DeclareMathOperator{\Imag}{Im}

\begin{document}

\title[Moment maps, symplectomorphisms and complex structures]
{Moment maps, symplectomorphism groups and compatible complex structures}
\author{Miguel Abreu}
\address{Departamento de Matem\'atica, Instituto Superior T\'ecnico,   
Av. Rovisco Pais,\newline\indent 1049-001 Lisboa, Portugal}
\email{mabreu@math.ist.utl.pt, ggranja@math.ist.utl.pt}
\author{Gustavo Granja}
\author{Nitu Kitchloo}
\address{Department of Mathematics, University of California, San
  Diego, USA}
\email{nitu@math.ucsd.edu}
\thanks{Nitu Kitchloo is supported in part by NSF through grant DMS
  0436600. Miguel Abreu and Gustavo Granja are supported in part by 
  FCT through program POCTI-Research Units Pluriannual Funding Program 
  and grants POCTI/MAT/57888/2004 and POCTI/MAT/58497/2004.}

\date{\today}

{\abstract In this paper we apply Donaldson's general moment map framework for the action
of a symplectomorphism group on the corresponding space of compatible (almost)
complex structures to the case of rational ruled surfaces. This gives a new
approach to understanding the topology of their symplectomorphism
groups, based on a result of independent interest: the space of compatible
integrable complex structures on any symplectic rational ruled surface is
(weakly) contractible. We also explain how in general, under this condition,
there is a direct relationship between the topology of a symplectomorphism
group, the deformation theory of compatible complex structures and the groups of
complex automorphisms of these complex structures.
}


\maketitle

\section{Introduction}

The known results regarding the topology of symplectomorphism groups,
for $2$-dimensional surfaces and $4$-dimensional rational ruled
surfaces, indicate a direct relation with the topology of the
corresponding groups of complex automorphisms. The goal of this paper
is to formulate this empirical relation more precisely, within
a general framework involving infinite dimensional groups, 
manifolds and moment maps. This general framework goes back to 
Atiyah and Bott~\cite{AB}, was made rigorous and precise in finite 
dimensions by Kirwan~\cite{K}, and was shown to apply to the action
of a symplectomorphism group on the corresponding space of compatible 
(almost) complex structures by Donaldson~\cite{D1}. 

It is usually quite difficult to make this general framework rigorous
and precise in infinite dimensions, and we will make no attempt at
that here. However, its formal application is often very useful as a
guide to the results one should expect and to the approach one might
take in proving them. We will present here two results,
Theorems~\ref{thm:contcs} and~\ref{thm:cohbsymp}, that are an example
of this. Their rigorous proofs will appear in~\cite{AGK}.

The paper is organized as follows. In section~\ref{s:moment} we
describe the general moment map framework and how it applies to the
action of a symplectomorphism group on the corresponding space of
compatible (almost) complex structures. In section~\ref{s:rrs} we
discuss the particular case of rational ruled surfaces,
where the geometric picture suggested by the moment map framework is
quite accurate and gives a new approach to understanding the
topology of their symplectomorphism groups. This approach is based on
the fact that the space of compatible integrable complex structures on
any symplectic rational ruled surface is (weakly) contractible
(Theorem~\ref{thm:contcs}). Any time this condition holds, there
should be a direct relationship between the problem of
understanding the topology of the symplectomorphism group and the
following two problems: understanding the deformation theory of
compatible complex structures and understanding the topology of the
groups of complex automorphisms of these complex structures. This
general relation is explained in section~\ref{s:vague}.

\subsection{Acknowledgements} Miguel Abreu and Gustavo Granja are
members of the Center for Mathematical Analysis, Geometry and
Dynamical Systems, which supported visits of Nitu Kitchloo to
IST-Lisbon. We thank the Center for its research environment and
hospitality. Miguel Abreu also thanks the organizers of the Conference
on Symplectic Topology at Stare Jablonki, Poland, for the opportunity
to speak about this work at such a wonderful event.

\section{Moment Map Geometry} \label{s:moment}

In this section, following~\cite{D1}, we recall a general moment map
framework and how it applies to the action of a symplectomorphism
group on the corresponding space of compatible almost-complex
structures.

\subsection{General Framework}

Let $G$ be a Lie group, $\Gg$ its Lie algebra,
$\langle\cdot,\cdot\rangle$ an inner product on $\Gg$ invariant under
the adjoint action, $\Gg^\ast$ the dual Lie algebra naturally
identified with $\Gg$ via $\langle\cdot,\cdot\rangle$, and $G^\C$ a
complexification of $G$. 

Let $(X,J,\Om)$ be a K\"ahler manifold equipped with an action of $G$
by K\"ahler isometries, i.e. a homomorphism
\[ \rho : G \to \Iso (X,J,\Om) = \Hol (X,J) \cap \Symp (X, \Om)\ .\]
Suppose this action satisfies the following two conditions:
\begin{enumerate}
\item[(i)] the holomorphic action of $G$ on $(X,J)$ extends to a
  holomorphic action of $G^\C$ on $(X,J)$;
\item[(ii)] the symplectic action of $G$ on $(X,\Om)$ admits a
  suitably normalized equivariant moment map $\mu : X \to \Gg^\ast$.
\end{enumerate}
Then we have the following two general principles.  

\vspace*{1ex}

\noindent {\bf General Principle I.} The complex and symplectic
quotients of $X$ by $G$ are naturally identified. More precisely,
\[ \mu^{-1}(0) / G \,=\, X^s / G^\C \]
where $X^s \subset X$ is an open, $G^\C$-invariant, subset of ``stable
points''. We will not define here this notion of stability, the
important point being that it should only depend on the holomorphic
geometry of the situation. The content of this principle is that 
on each stable $G^\C$-orbit there is a point $p\in\mu^{-1}(0)$, 
unique up to the action of $G$. 

\vspace*{1ex}

\noindent {\bf General Principle II.} The map $\|\mu\|^2 \equiv
\langle \mu,\mu \rangle : X \to \R$ behaves like a $G$-invariant
Morse-Bott function, whose critical manifolds compute the equivariant
cohomology $H^\ast_G (X) \equiv  H^\ast (X \times_G EG)$
(over $\Q$ and, in good special cases, also over $\Z$).

\vspace*{1ex}

Combining these two general principles, one gets the following
geometric picture for the action of $G$ on $X$:
\begin{enumerate}
\item[-] The gradient flow of $-\|\mu\|^2$ induces an invariant
  stratification 
  \[ X = V_0 \sqcup V_1 \sqcup V_2 \sqcup \cdots \ ,\]
  where each $V_k$ is the stable manifold of some critical set
  $C_k$ of $\|\mu\|^2$.
\item[-] Let $\Oo_k$ denote the coadjoint orbit $G\cdot \xi_k \subset
  \Gg^\ast$, where $\xi_k = \mu (p_k)$ for some $p_k \in C_k$. Then
  \[ V_k / G^\C \, \simeq \, \mu^{-1}(\Oo_k)/G.\ \]
  If $C_0 = \mu^{-1}(0)$ then $\Oo_0 = \{ 0 \}$ and $V_0 = X^s$.
\item[-] The equivariant cohomology $H^\ast_G (X)$ can be computed
  from $H^\ast_G (V_k)\,,\ k=0,1,2,\ldots$ 
  (over $\Q$ and, in good special cases, also over $\Z$).
\end{enumerate}

\subsection{Symplectomorphism Groups and Compatible Complex Structures}

Consider a compact symplectic manifold $(M,\om)$, of dimension $2n$,
and assume that $H^1(M,\R) = 0$. Let $G \equiv \Symp (M,\om)$ be the
symplectomorphism group of $(M,\om)$. This is an infinite dimensional Lie
group whose Lie algebra $\Gg$ can be identified with the space of
functions on $M$ with integral zero:
\[ \Gg = C^\infty_0 (M) \equiv \left\{ f:M\to\R\,:\ \int_M f\ 
  \frac{\om^n}{n!} \, = \, 0 \right\}\ .\]
$\Gg$ has a natural invariant inner product
$\langle\cdot,\cdot\rangle$, given by
\[ \langle f,g\rangle \equiv \int_M f\cdot g\  \frac{\om^n}{n!}\ ,\]
which will be used to identify $\Gg^\ast$ with $\Gg$.

Consider now the space $\Jj (M,\om)$ of almost complex structures $J$
on $M$ which are compatible with $\om$, i.e. for which the bilinear
form
\[ g_{J} (\cdot,\cdot) = \om (\cdot,J\cdot)\]
is a Riemannian metric on $M$. This is the space of sections of a
bundle over $M$ with fiber the contractible symmetric K\"ahler
manifold $Sp(2n,\R)/U(n) \equiv$ Siegel upper half
space~\cite{S}. This fiberwise symmetric K\"ahler structure, together
with the volume form induced by $\om$ on $M$, turns $\Jj(M,\om)$ into
an infinite dimensional (contractible) K\"ahler manifold.

The symplectomorphism group $G$ acts naturally on $\Jj(M,\om)$ by
K\"ahler isometries:
\[ \phi \cdot J \equiv \phi_\ast (J) = d\phi \circ J \circ
d\phi^{-1}\,,\ \forall \phi\in G\,,\ J\in\Jj(M,\om)\ .\]
To fit the previous general framework, this action should satisfy
conditions (i) and (ii).

The first (holomorphic) condition poses an immediate problem since
there is no complexification $G^\C$ of the symplectomorphism group
$G$. However, we can certainly complexify the Lie algebra $\Gg$ to
\[ \Gg^\C \equiv \left\{ f:M\to\C\,:\ \int_M f\ 
  \frac{\om^n}{n!} \, = \, 0 \right\}\]
and the infinitesimal action of $\Gg$ on $\Jj(M,\om)$ extends to an
action of $\Gg^\C$, since the complex structure on $\Jj(M,\om)$ is
integrable. This gives rise to an integrable complex distribution on
$\Jj(M,\om)$ whose leaves play the role of ``connected components of 
orbits of the group $G^\C$''. 

In the holomorphic side of General Principles I and II that we want to
apply, $G^\C$ is not that important when compared with the role played
by its orbits. The geometric meaning of these ``$G^\C$-orbits''
becomes quite clear if one restricts the actions under consideration
to the invariant K\"ahler submanifold $X$ of compatible integrable
complex structures
\[ X \equiv \Jj^{\rm int} (M,\om) \subset \Jj (M,\om)\,,\]
determined by the vanishing of the Nijenhuis tensor. Here, it follows
from Donaldon's analysis in~\cite{D1} that

  $J,J'\in X$ belong to the same ``$G^\C$-orbit'' iff there
  exists $\varphi\in\Diff (M)$ such that 
  \[ [\varphi^\ast (\om)] = [\om] \in H^2 (M,\R) \quad\text{and}\quad
  \varphi_\ast (J) = J'\ .\]
This explicit description of a ``$G^\C$-orbit'' is good enough to
consider that the holomorphic action of $G = \Symp (M, \om)$ on $X =
\Jj^{\rm int} (M,\om)$ satisfies condition (i).

Regarding condition (ii), Donaldson~\cite{D1} shows that there always
exists an equivariant and suitably normalized moment map
\[ \mu : \Jj(M,\om) \to \Gg^\ast \cong C^\infty_0 (M) \]
for the symplectic action of $G$ on $\Jj(M,\om)$, given by
\[ \mu (J) = \left(\text{Hermitian scalar curvature $S(J)$ of the metric
  $g_J$}\right) - d\,,\]
where $d$ is the constant defined by
\[
d\cdot \int_M \frac{\om^n}{n!} \equiv 2\pi c_1(M) \wedge
[\om]^{n-1} (M) = \int_M S(J) \ \frac{\om^n}{n!}\,.
\]
Note that on $X\subset \Jj (M,\om)$, i.e. for integrable $J$, the Hermitian
scalar curvature $S(J)$ coincides with the usual scalar curvature of
the Riemannian metric $g_J$. 

We have concluded that the K\"ahler action of $G$ on $X$ satisfies
conditions (i) and (ii) of the general framework, and so General
Principles I and II should apply. What do they say in this context?

\vspace*{1ex}

\noindent {\bf General Principle I} says that each stable compatible
complex structure is diffeomorphic to one in $\mu^{-1}(0)$, unique up
to the action of $G$. Since
\[ J\in\mu^{-1}(0) \Leftrightarrow \mu(J) = 0 \Leftrightarrow
S(J) = d = \,\text{constant,}\]
this says that on each diffeomorphism class of stable compatible complex
structures there should exist a unique $\Symp(M,\om)$-orbit
whose corresponding K\"ahler metric has constant scalar curvature.
(See the work of Donaldson~\cite{D2}, \cite{D3} and
\cite{D4},  exploring this consequence of General Principle I.)

\vspace*{1ex}

\noindent {\bf General Principle II} says that the critical points of
\[ \|\mu\|^2 : X = \Jj^{\rm int} (M, \om) \to \R\,,\ 
\|\mu\|^2 (J) = \int_M S^2(J)\ \frac{\om^n}{n!}\,+\,\text{constant,}\]
determine the equivariant cohomology $H^\ast_G (X)$. 

These critical points are, in particular, extremal K\"ahler metrics in the 
sense of Calabi (\cite{C1} and~\cite{C2}). When extremal K\"ahler metrics exist, they minimize
$\|\mu\|^2$ on the corresponding "$G^\C$-orbit" (see~\cite{H}) and are conjecturally
unique up to the action of $G$ (see~\cite{CC}).

\vspace*{1ex}

In the concrete examples we will discuss (Riemann surfaces and rational
ruled surfaces) these general principles do hold. Whenever this is the case,
one gets the following geometric picture for the action of $G = \Symp
(M,\om)$ on $X=\Jj^{\rm int} (M,\om)$:
\begin{enumerate}
\item[-] There is a stratification of $X$ of the form
\[ X = V_0 \sqcup V_1 \sqcup V_2 \sqcup \cdots\,.\]
In this stratification, each $V_k$ contains the set of
compatible complex structures which are diffeomorphic to an extremal
one with normalized scalar curvature in
\[\Oo_k \equiv\,\text{coadjoint orbit $G\cdot (S(J_k)-d)$ in $C^\infty_0
  (M) \cong \Gg^\ast$,}\]
where $J_k$ is some critical compatible complex structure in
$V_k$. 
\item[-] ``$V_k / G^\C$''$\equiv \{$ ``$G^\C$-orbits'' in $V_k \} \simeq
  \mu^{-1} (\Oo_k) / G$ is some moduli space of complex
  structures, that one might try to understand using methods from
  complex geometry (deformation theory).
\item[-] Let $\Oo_{J_k}$ denote the ``$G^\C$-orbit'' through some
  extremal $J_k \in V_k$ and let $K_k \equiv \Iso (M,\om,J_k) \subset
  G$. Then, if the group $\Hol_{[\om]}(M,J_k)$ of holomorphic automorphisms
  which preserve the cohomology class of $\om$ is the complexification of $K_k$
  (this is always the case if the groups are connected by~\cite{C2}) then the inclusion
  \[G/K_k \cong G\cdot J_k \hookrightarrow \Oo_{J_k} = 
  \text{``$G^\C$''} / \Hol_{[\om]} (M, J_k)\]
  is a homotopy equivalence.
\item[-] The equivariant cohomology $H^\ast_G (X)$ can be computed
  from $H^\ast_G (V_k)\,,\ k=0,1,2,\ldots$. From the previous two
  points, each $H^\ast_G (V_k)$ should be determined from finite
  dimensional considerations involving moduli spaces of complex
  structures and subgroups of isometries in $G = \Symp (M,
  \om)$. Recall that $\Jj (M,\om)$ is always contractible. If $X =
  \Jj^{\rm int} (M,\om) \subset \Jj (M,\om)$ is also contractible,
  then
  \[ H^\ast_G (X) = H^\ast (BG) \ .\]
\end{enumerate}

\section{Rational Ruled Surfaces} \label{s:rrs}

In this section we discuss the particular case of rational ruled
surfaces, formulating the precise results suggested by
the framework of section~\ref{s:moment}.

\subsection{Symplectic Structures}

As smooth $4$-manifolds, rational ruled surfaces are $S^2$-bundles
over $S^2$. Since $\pi_2 \left(B\Diff^+(S^2)\right) \cong 
\pi_2 \left(BSO(3)\right) \cong \Z/2$, there are only two
diffeomorphism types classified by the second Stiefel-Whitney class of
the bundle (the mod $2$ reduction of the Euler class): a trivial
$S^2$-bundle over $S^2$ and a nontrivial $S^2$-bundle over
$S^2$. Since the story for each of these is analogous, we will
concentrate here on the trivial bundle, i.e.
\[ M = S^2 \times S^2\,.\]
Symplectic structures on $S^2 \times S^2$ are classified by the
following theorem.
\begin{thm} \label{thm:lm} \emph{(Lalonde-McDuff)}
If $\om$ is a symplectic form on $S^2 \times S^2$, then it is
diffeomorphic to $\la \sigma \oplus \mu \sigma$ for some real $\la,\mu
> 0$, where $\sigma$ denotes the standard area form on $S^2$ with
$\int_{S^2}\,\sigma = 1$.
\end{thm}
Since the symplectomorphism group and the space of compatible almost
complex structures is not affected by positive scalings of the
symplectic form and we can switch the two $S^2$-factors, it will
suffice to consider symplectic structures of the form
\[
\om_\la = \la \sigma \oplus \sigma \quad\text{with}\quad
1\leq \la \in\R\,.
\]
From now on we will use the following notation:
\begin{align}
M_\la & = (S^2\times S^2, \om_\la)\,,\ 1\leq \la\in\R\,; \notag \\
G_\la & = \Symp (M_\la) = \text{symplectomorphisms of $M_\la$;} \notag \\
\Jj_\la & = \Jj (M_\la) = \text{compatible almost complex structures;}
\notag \\
X_\la & = \Jj^{\rm int} (M_\la) = \text{compatible integrable 
complex structures.} \notag
\end{align}
We will also use the following obvious isomorphism:
\begin{align}
H_2 (S^2\times S^2, \Z) & \overset{\cong}{\lra} \Z \oplus \Z \notag \\
m \left[S^2\times \text{pt}\right] + n \left[\text{pt}\times S^2\right]
& \longmapsto (m,n) \ . \notag
\end{align}

\subsection{Compatible Integrable Complex Structures}

As a complex manifold, a rational ruled surface is a holomorphic
$\CP^1$-bundle over $\CP^1$. These are the well known Hirzebruch
surfaces
\[
H_k = P (\Oo\oplus\Oo(-k)) \quad\text{for some $k\in\N_0$,}
\]
where we write $\Oo(-1)$ for the tautological line bundle over $\CP^1$
and $P(E)$ for the projectivization of a vector bundle $E$.

Any complex structure $J$ on $S^2 \times S^2$ is isomorphic to
$H_{2k}$ for some $k\in\N_0$, while the ``odd'' Hirzebruch surfaces
are diffeomorphic to the nontrivial $S^2$-bundle over $S^2$
(see~\cite{Q}). When $(S^2\times S^2, J)$ has two embedded $\CP^1$'s
with self-intersection $0$ and themselves intersecting at one point,
then $(S^2\times S^2, J) \cong H_0$. When $(S^2\times S^2, J)$ has an
embedded $\CP^1$ with self-intersection $-2k < 0$, then 
$(S^2\times S^2, J) \cong H_{2k}$.

To understand which of these complex structures $J$ can be made
compatible with a symplectic form $\om_\la$, for some $1\leq\la\in\R$,
it is important to note that the compatibility condition implies that
the symplectic form evaluates positively on any $J$-holomorphic
curve. Hence, for a compatible $J\in X_\la$, a homology class
$(m,n)\in H_2 (S^2\times S^2;\Z)$ can only be represented by a
$J$-holomorphic curve if $\la m + n > 0$. This rules out embedded curves with
self-intersection less than $-2\ell$, where $\ell\in\N_0$ is such that
$\ell < \la \leq \ell +1$. In particular, the class $(1,-k)\in H_2
(S^2\times S^2; \Z)$, with self-intersection $-2k$, can only be
represented by a $J$-holomorphic curve for some $J\in X_\la$ if $\la -
k > 0$.

This turns out to be the only relevant condition. In fact, we have the
following theorem.
\begin{thm} \label{thm:cs}
Given $1\leq \la \in \R$, there is a stratification of $X_\la$ of the
form
\[ X_\la = V_0 \sqcup V_1 \sqcup \cdots \sqcup V_\ell\,,\]
with $\ell\in\N_0$ such that $\ell < \la \leq \ell +1$ and where:
\begin{enumerate}
\item[(i)]
\begin{align}
V_k \equiv \{J\in X_\la &:\,(S^2\times S^2) \cong H_{2k} \} 
\notag \\
= \{J\in X_\la &:\,\text{$(1,-k)\in H_2 (M,\Z)$ is
    represented} \notag \\
& \quad\text{by a $J$-holomorphic sphere}\}\,. \notag
\end{align}
\item[(ii)] $V_0$ is open and dense in $X_\la$. For $k\geq 1$, $V_k$
  has codimension $4k-2$ in $X_\la$.
\item[(iii)] $\overline{V_k} = V_k \sqcup V_{k+1} \sqcup \cdots \sqcup V_\ell$.
\item[(iv)] For each $k\in\N_0$, there is a complex structure $J_k\in
  V_k$, unique up to the action of $G_\la$, for which $g_{\la,k} \equiv \om_\la
  (\cdot, J_k \cdot )$ is an extremal K\"ahler metric.
\item[(v)] Denoting by $K_k$ the K\"ahler isometry group of
  $(S^2\times S^2, \om_\la, J_k)$, we have that
\[
K_k \cong 
\begin{cases}
\Z/2 \ltimes (SO(3)\times SO(3))\,, &\text{if $k=0$;} \\
S^1\times SO(3)\,, &\text{if $k\geq 1$.}
\end{cases} 
\]
\item[(vi)] Given $J\in V_k$, there exists $\varphi \in \Diff (S^2
  \times S^2)$ such that 
\[ [\varphi^\ast (\om_\la)] = [\om_\la] \in H^2 (S^2\times S^2;\R) 
\quad\text{and}\quad \varphi_\ast (J_k) = J\]
so each stratum $V_k$ consists of a unique ``$G_\la^\C$''-orbit, the
orbit through $J_k$. Moreover, the inclusion
\begin{align}
\left(G_\la / K_k \right) & \lra V_k = \Oo_{J_k} \notag
\\
[\psi] & \longmapsto \psi_\ast (J_k) \notag
\end{align}
is a weak homotopy equivalence.
\item[(vii)] Each $V_k$ has a tubular neighborhood $NV_k \subset
  X_\la$, with normal slice given by
\[H^1 (H_{2k}, \Theta) \cong \C^{2k-1}\,,\]
where $\Theta =$ sheaf of holomorphic vector fields on $H_{2k}$.
\item[(viii)] For $k\geq 1$, the representation of $K_k \cong S^1
  \times SO(3)$ on the normal slice $\C^{2k-1}$ at $J_k \in V_k$ is the following:
  $S^1$ acts diagonally and $SO(3)$ acts irreducibly with highest
  weight $2(k-1)$.
\end{enumerate}
\end{thm}
\begin{proof}
Points (i), (ii), (iii), (v) and (vii) follow from standard complex
geometry and deformation theory applied to complex structures on $S^2
\times S^2$ (see~\cite[I.6]{Ca}). 
One needs to check that standard deformation theory can
in fact be used here, in the context of compatible complex
structures. This, together with points (vi) and (viii), will be proved
in~\cite{AGK}. Point (iv) is proved in~\cite{C1}.
\end{proof}

This theorem shows that the geometric picture suggested by the moment
map framework of section~\ref{s:moment} is quite accurate for rational
ruled surfaces.  It implies, by standard equivariant cohomology theory
(see~\cite{AGK}), the following corollary.
\begin{cor} \label{cor:cs}
Given $\ell\in\N_0$ and $\la\in\left]\ell,\ell+1\right]$, we have that
\[
H^\ast_{G_\la} (X_\la;\Z) \cong H^\ast (BSO(3)\times BSO(3);\Z) \oplus
\oplus_{k=1}^\ell \Sigma^{4k-2} H^\ast (BS^1 \times BSO(3);\Z)\,,
\]
where $\cong$ indicates a group isomorphism.
\end{cor}

\subsection{Compatible Almost Complex Structures}

We will now look at the contractible space $\Jj_\la$ of compatible
almost complex structures on $M_\la = (S^2 \times S^2, \om_\la)$, for
a given $1\leq \la \in \R$. Gromov~\cite{G} was the first to study
the topology of the symplectomorphism group $G_\la$ by analysing its
natural action on this space $\Jj_\la$. He used pseudo-holomorphic
curves techniques and this approach turned out to be quite fruitful
(see~\cite{A} and~\cite{AM}). In fact, the action of $G_\la$
is compatible with a natural geometric stratification of $\Jj_\la$,
analogous to the one presented in Theorem~\ref{thm:cs} for $X_\la
\subset \Jj_\la$. Replacing holomorphic spheres by
pseudo-holomorphic ones, and complex deformation theory by gluing
techniques for pseudo-holomorphic spheres, one can prove the
following theorem (see~\cite{M}).
\begin{thm} \label{thm:acs}
Given $1\leq \la \in \R$, there is a stratification of the
contractible space $\Jj_\la$ of the form
\[ \Jj_\la = U_0 \sqcup U_1 \sqcup \cdots \sqcup U_\ell\,,\]
with $\ell\in\N_0$ such that $\ell < \la \leq \ell +1$ and where:
\begin{enumerate}
\item[(i)]
\begin{align}
U_k \equiv \{J\in \Jj_\la &:\,\text{ $(1,-k)\in H_2 (M;\Z)$ is
    represented} \notag \\
& \quad\text{by a $J$-holomorphic sphere}\}\,. \notag
\end{align}
In particular, $V_k = U_k \cap X_\la \subset U_k$.
\item[(ii)] $U_0$ is open and dense in $\Jj_\la$. For $k\geq 1$, $U_k$
  has codimension $4k-2$ in $\Jj_\la$.
\item[(iii)] $\overline{U_k} = U_k \sqcup U_{k+1} \sqcup \cdots \sqcup U_\ell$.
\item[(iv)] The inclusion
\begin{align}
\left(G_\la / K_k \right) & \lra U_k \notag
\\
[\psi] & \longmapsto \psi_\ast (J_k) \notag
\end{align}
is a weak homotopy equivalence, where $J_k \in V_k \subset U_k$ and
$K_k = \Iso (S^2\times S^2, \om_\la, J_k)$ were characterized in
Theorem~\ref{thm:cs}. In particular, the inclusion 
\[
V_k \lra U_k
\]
is also a weak homotopy equivalence.
\item[(v)] Each $U_k$ has a tubular neighborhood $NU_k \subset
  \Jj_\la$ which fibers over $U_k$ as a ball bundle.
\end{enumerate}
\end{thm}

\subsection{Contractibility of $X_\la$}

Given $\ell \in\N_0$ and $1\leq \la \in \left]\ell, \ell+1\right]$, we
can combine the results of Theorems~\ref{thm:cs} and~\ref{thm:acs} to
obtain a finite family of diagrams, one for each 
$0\leq k \leq \ell$, of the form
\[
\xymatrix{  
F_k \ar[r] \ar[d]  & NV_k \ar[r] \ar[d] & 
V_k \subset X_\la \ar[d] \\
F_k \ar[r]  & NU_k \ar[r] & U_k \subset \Jj_\la    
}
\]
where the vertical arrows are inclusions, the one on the left 
representing the identity between the fibers of the tubular 
neighborhoods over $V_k \subset U_k$. 
These diagrams are $G_\la$-equivariant in a suitable sense.
Given that $\Jj_\la$ is
contractible and $V_k$ is weakly homotopy equivalent to $U_k$, one 
can use this finite family of diagrams to prove the 
following theorem (see~\cite{AGK}).
\begin{thm} \label{thm:contcs}
Given $1\leq\la\in\R$, the space $X_\la$ of compatible integrable
complex structures on $(S^2 \times S^2, \om_\la)$ is weakly
contractible. 
\end{thm}
As far as we know, these are the first known examples of dimension
greater than two where the topology of the space of compatible
integrable complex structures is understood.

\subsection{Cohomology of $BG_\la$}

Theorem~\ref{thm:contcs} implies that
\[ H^\ast_{G_\la} (X_\la;\Z) \cong H^\ast (BG_\la; \Z)\,.\]
Combining this isomorphism with Corollary~\ref{cor:cs}, we get the
following theorem.
\begin{thm} \label{thm:cohbsymp}
Given $\ell\in\N_0$ and $\la\in\left]\ell,\ell+1\right]$, we have that
\[
H^\ast (BG_\la;\Z) \cong H^\ast (BSO(3)\times BSO(3);\Z) \oplus
\oplus_{k=1}^\ell \Sigma^{4k-2} H^\ast (BS^1 \times BSO(3);\Z)\,,
\]
where $\cong$ indicates a group isomorphism.
\end{thm}
Although we used a new point of view, this theorem
is not the first result regarding the topology of $BG_\la$:
\begin{itemize}
\item[-] the rational cohomology ring of $BG_\la$ was determined 
  in~\cite{AM} for any $1\leq\la\in\R$;
\item[-] the integral cohomology ring of $BG_\la$ was determined for
  $\la = 1$ in~\cite{G} (where Gromov proves that $G_1$ is homotopy 
  equivalent to $SO(3)\times SO(3)$) and for $1<\la\leq 2$
  in~\cite{AG}. 
\end{itemize}
We hope to include in~\cite{AGK} a description of the ring structure
of $H^\ast (BG_\la;\Z)$ for any $1\leq\la\in\R$.

\section{Symplectomorphisms and complex automorphisms} \label{s:vague}

Let $(M,\omega)$ be a symplectic manifold and $\Jj^{\rm int}(M,\omega)$
the corresponding space of compatible integrable complex structures.
The aim of this section is to explain how the condition 
\begin{equation}
\label{contract}
\Jj^{\rm int}(M,\omega) \sim \ast
\end{equation} 
that the space of compatible complex structures be contractible,
essentially equates the problem of understanding the topology of 
$\Symp(M,\omega)$ with the following two problems
\begin{itemize}
\item Understanding the (large-scale) deformation theory of complex 
structures compatible with $\omega$.
\item Understanding the topology of the groups of complex
  automorphisms of each of these complex structures.
\end{itemize}
The same is true if one replaces compatible complex structures by 
tame complex structures in which case the deformation theory in
question is the usual deformation theory of Kodaira and Spencer 
(as the tameness condition is open).

We should point out straight away that we only know condition 
\eqref{contract} to hold in a few very simple examples: 
Riemann surfaces (where the condition holds because all almost complex 
structures are integrable) and simple rational complex surfaces such
as $\CP^2$ and rational ruled surfaces. 
The work of Lalonde and Pinsonnault \cite{LP} suggests that condition 
\eqref{contract} also holds for blow-ups of $\CP^2$ at two points. It 
would be very interesting to have some understanding of the generality
of \eqref{contract}.

The reason why \eqref{contract} equates the problems above is that it 
establishes a sort of weak equivalence between two moduli problems - 
one in symplectic geometry and another in complex geometry. 
To explain this we must first recall some basic facts regarding 
topological groupoids. The reader is referred to \cite[Section I]{Ha} for more
explanation and details.

\subsection{Topological groupoids}

A topological groupoid is a small topological category in which every 
morphism is invertible. Thus, a topological groupoid consists of a
pair of spaces $(O,M)$ together with some structure maps. $O$ is the
space of objects and $M$ the space of isomorphisms between the objects
in $O$. The structure maps are $\iota: O \to M$ assigning the identity 
isomorphism to each object, two maps $d_0,d_1:M \to O$ assigning to a 
morphism its domain and range, a composition map 
$\mu: M\times_O M \to M$ and an inverse map $c: M \to M$ 
satisfying the obvious identities.

To a groupoid one can associate the space of isomorphism 
classes $O/I$, which we will call the \emph{coarse moduli space},
 defined as the quotient of $O$ by the equivalence relation
generated by $d_0(\alpha)\sim d_1(\alpha)$. 
Of course the groupoid carries much more information than the quotient space.

In many cases, including all those that will concern us, 
moduli problems can be described in terms of topological
groupoids of the following special form:
\begin{example}
\label{action}
Let $G$ be a topological group acting on a space $X$. We define a
groupoid $\Gamma(X,G)$ with objects 
$O=X$ and $M=G\times X$. The structure maps are 
$\iota(x)=(1,x)$, $d_0(g,x) = x$, $d_1(g,x) = g \cdot x$,
$\mu(g,x,h,g\cdot x) = (hg,x)$ and $c(g,x) = (g^{-1},g\cdot x)$. 
The associated quotient space is the space
of orbits $X/G$.

Note that $\Gamma(\ast,G)$ can be naturally identified with the topological group $G$.
\end{example}

Let $\Gamma=(O,M)$ be a topological groupoid, $B$ be a space and 
$\cat U=\{U_i\}$ be an open cover of $B$.
A 1-cocycle over $\cat U$ with values in $\Gamma$ is the assignment 
to each pair $i,j$ of a map
\[  g_{ij} : U_i \cap U_j \to M \]
satisfying the cocycle condition $\mu (g_{ij} , g_{jk}) = g_{ik}$ 
on triple intersections. In particular, $g_{ii}: U_i \to M$ sends each 
point of $U_i$ to an identity morphism and hence amounts to an object 
selection map $f_i: U_i \to O$. The cocycle condition then says that 
these selections are compatible so that the
$1$-cocycle represents a continuous family of objects over $B$. 
One defines equivalence between $1$-cocycles in the obvious way. 
An equivalence class of $1$-cocycles is called a 
\emph{$\Gamma$-structure on $X$}. 

Two $\Gamma$-structures $\sigma_i$, $i=0,1$, on $Y$ are said to be 
\emph{concordant}\footnote{The term used by Haefliger is homotopic.} 
if there is a $\Gamma$-structure $\sigma$ on $Y\times[0,1]$ such that 
$\sigma_{|Y\times 0}=\sigma_0$ and $\sigma_{|Y\times 1} = \sigma_1$.  
If $\Gamma=\Gamma(\ast,G)$ is a topological group, concordance is 
the same as equality but this is generally not the case.

Any topological groupoid $\Gamma$ has a classifying space 
$B\Gamma$ which is a homotopy invariant 
version of the coarse moduli space. It is determined by the equation
\[ [Y,B\Gamma] = \{ \text{ concordance classes of } 
\Gamma-\text{structures on } Y \} \]
In the case of Example \ref{action} the classifying space of
$\Gamma(X,G)$ is the 
Borel construction on the $G$-space $X$. 
That is, if we write $EG$ for a contractible free $G$-space, we have
\[ B\Gamma(X,G) \simeq EG\times_G X \]
In particular, if the $G$-space $X$ consists of a single orbit and $H$ 
denotes the isotropy subgroup of one of the points of $X$, we have 
\begin{equation}\label{bh} 
B\Gamma \simeq EG \times_G G/H = EG/H \simeq BH 
\end{equation}
and so the classifying space of the groupoid is the same as the
classifying space of the isotropy group
(of any object). If there is more than one orbit then the map 
$EG \to \ast$ still induces a map
\[ \pi: B\Gamma \to X/G \]
with fibers
\[ \pi^{-1}(Gx) \simeq B\Aut(x) \]
where we have written $\Aut(x)$ for the isotropy group of $x$.
Intuitively, this says that the space $B\Gamma$ is obtained by
 gluing the classifying spaces of the automorphism
groups $B\Aut(x)$ via the topology of the moduli space $X/G$, 
which one can write 
\begin{equation}
\label{bgamma}
B\Gamma \simeq \int_{X/G} B\Aut(x). 
\end{equation}
In good situations, this is a precise statement \cite[Proposition, p. 183]{DF}. 
Even if it is not, it still provides a 
useful guide to understanding the topology of $B\Gamma$. 
We will see how this plays out in Examples~\ref{ell} and~\ref{rr} below.

We will say that a map $\phi: \Gamma_1 \to \Gamma_2$ of 
topological groupoids is a \emph{weak equivalence}
if the induced map of classifying spaces
\[ B\phi : B\Gamma_1 \to B\Gamma_2 \]
is a weak homotopy equivalence. In the examples we are concerned with, 
the source of weak equivalences is
the following immediate consequence of the homotopy invariance of the 
Borel construction:
\begin{lemma} \label{morita}
If $X$ and $Y$ are $G$-spaces and $f:X \to Y$ is a $G$-equivariant 
map which is a weak equivalence then the induced map of groupoids 
$\Gamma f: \Gamma(X,G) \to \Gamma(Y,G)$ is a weak equivalence.
\end{lemma}

\subsection{Examples of moduli problems}

Let $M$ be a compact manifold and $G=\Diff^+(M)$ the group of 
orientation preserving diffeomorphisms of $M$ with the $C^\infty$ topology.
Consider the following $G$-spaces:
\begin{enumerate}[(i)]
\item The space $\Jj^{\rm int}(M)$ of complex structures on $M$
  compatible with some (unspecified) symplectic form on $M$.
\item The space $\Omega^{int}(M)$ of symplectic forms on $M$.
\item The space $K^{int}(M) =\{ (J,\omega) \in \Jj^{\rm int}(M) 
\times \Omega^{int}(M) \colon J \text{ compatible with } \omega \}$ 
of K\"ahler structures on $M$.
\end{enumerate}
We will write $\Gamma \Jj$, $\Gamma \Omega$ and $\Gamma K$ 
for the associated topological groupoids.
The projections 
\[ 
\Jj^{\rm int}(M) 
\llla{\pi_1} K^{int}(M) \llra{\pi_2} \Omega^{int}(M) 
\]
are $G$-equivariant. The fibers of $\pi_1$ are always convex so 
$\pi_1$ is a weak equivalence. Lemma \ref{morita} then says that 
the groupoids $\Gamma \Jj$ and $\Gamma K$ are weakly equivalent.
Condition $\eqref{contract}$ is precisely the condition that the
fibers of $\pi_2$ be contractible. If in addition
we know that $\pi_2$ is surjective, i.e. that all symplectic
structures on $M$ are K\"ahler, then Lemma \ref{morita}
says that $\Gamma \Omega$ is also weakly equivalent to 
$\Gamma \Jj$ and $\Gamma K$.

\begin{remark}
We can also consider the non-integrable moduli problems derived from
the space $\Jj(M)$ of almost complex structures compatible with some 
symplectic form on $M$ and $K(M)$ of almost-K\"ahler structures.
These groupoids are always weakly equivalent to $\Gamma\Omega$. 
In this light, condition \eqref{contract} can be seen as the
requirement that the inclusions $\Jj^{\rm int}(M) \to \Jj(M)$ or 
$K^{int}(M) \to K(M)$ induce weak equivalences of topological groupoids.
\end{remark}

Now consider the following variation on the previous situation. 
Fix a symplectic form $\omega$ on $M$ and consider the orbit
\[ 
X_0 = \Diff^+(M) \cdot \omega \subset \Omega^{int}(M) \,.
\]
Let $X_1=\pi_1^{-1}(X_0) \subset K^{int}(M)$ and 
$X_2=\pi_2(X_1) \subset \Jj^{\rm int}(M)$ so that we have 
$G$-equivariant maps
\[ 
X_0 \llla{\pi_1} X_1 \llra{\pi_2} X_2 
\]
If \eqref{contract} holds for the symplectic form $\omega$ 
then as before all three groupoids are weakly
equivalent and using \eqref{bh} and \eqref{bgamma} we have
\begin{equation}
\label{mainpoint}
B\Symp(M,\omega) \simeq B\Gamma(X_0,G) \simeq B\Gamma(X_2,G) \simeq \int_{[J] \in Z} B\Aut(J) 
\end{equation}
where $Z=X_2/\Diff^+(M)$ is the coarse moduli space of complex
structures compatible with $\omega$, and $\Aut(J)$ is the group of 
biholomorphisms of the complex manifold $(M,J)$.

\begin{example} {\bf Elliptic curves.} \label{ell}
Let $M=S^1\times S^1$ and let $\omega$ be the standard symplectic
form. The topology of the diffeomorphism
groups of surfaces is well understood. The inclusion of the affine
orientation preserving automorphisms of $S^1\times S^1$ in 
$\Symp(M,\omega)$ is a weak equivalence, i.e.
\[    
SL_2(\Z) \ltimes (S^1 \times S^1) \simeq \Symp(M,\omega)\,.
\]
The complex structures on a torus (which are all compatible with 
any given area form) are parametrized
by points in the upper half plane up to the action of $SL(2,\Z)$ 
by fractional linear transformations.
Thus the space $Z$ in \eqref{mainpoint} is obtained from 
\[ \Omega = \{ z \in \C \colon -\tfrac 1 2 \leq \Real(z) 
\leq \tfrac 1 2, |z|\geq 1, \Imag(z)>0 \} \]
by identifying the boundary points through the maps $z \mapsto z+1$ and $z \mapsto -1/z$. 
If $E_z$ denotes the elliptic curve parametrized by $z \in \Omega$, we have 
\[ \Aut_\C(E_z) = \Aut(<1,z>) \ltimes E_z \]
where $\Aut(<1,z>)$  denotes the automorphisms of the lattice 
generated by $\{1,z\}$. Most lattices have
multiplication by $-1$ as their only non-trivial automorphism. 
The exceptions are $<1,i>$ which has $\Z/4$ as 
automorphism group and $<1, e^{\pi i/3} >$ which has $\Z/6$ as 
automorphism group.

Thus, in this case \eqref{mainpoint} says that $B\Symp(M,\omega)$ 
can be expressed as the double mapping cylinder (or homotopy pushout) of the maps 
\[ 
B( \Z/4 \ltimes (S^1\times S^1)) \leftarrow B(\Z/2 \ltimes (S^1\times S^1)) 
\rightarrow B(\Z/6 \ltimes (S^1\times S^1)) 
\] 
which amounts to the familiar decomposition of $SL(2,\Z)$ as the 
amalgam of $\Z/4$ and $\Z/6$ over $\Z/2$ \cite[1.5.3]{Se}.
\end{example}

\begin{example} {\bf $S^2\times S^2$.} \label{rr}
Let us resume section~\ref{s:rrs} from the point of view of this section.
A theorem of Qin~\cite{Q} says that any complex structure on $S^2\times
S^2$ is isomorphic to one of the even
Hirzebruch surfaces $H_{2k}$, with $k\geq 0$. Theorem~\ref{thm:lm} 
(Lalonde-McDuff) says that any symplectic form on $S^2\times S^2$ is,
up to scale, diffeomorphic to 
$\omega_\lambda = \lambda \sigma \oplus \sigma$ for some $\lambda \geq 1$.
The complex structures which are compatible with $\omega_\lambda$ are
the $H_{2k}$ with $k<\lambda$. So in this example the coarse moduli
space $Z$ is a finite set of points (with a non-Hausdorff topology)
and \eqref{mainpoint} expresses $B\Symp(M,\omega_\lambda)$ as a union 
of spaces with the homotopy type of
$B\Aut_\C(H_{2k})$ for the allowable values of $k$.

The way the different spaces $B\Aut_\C(H_{2k})$ fit together is 
controlled by the deformation theory of the complex structures, 
which in this case is the usual Kodaira-Spencer theory.
\end{example}

\end{document}